\documentclass[11pt,a4paper]{article}
\usepackage{latexsym,amsfonts}

\begin{document}

\newtheorem{theo}{Theorem}
\newtheorem{prop}[theo]{Proposition}
\newtheorem{lem}[theo]{Lemma}
\newtheorem{efi}[theo]{Definition}
\newtheorem{coro}[theo]{Corollary}

\newcommand{\R}{\ensuremath{\mathbb R}}

\newcommand{\qed}{\hfill $\Box$}

\newcommand{\beq}{\begin{eqnarray*}}
\newcommand{\eeq}{\end{eqnarray*}}

\newcommand{\lag}{\langle}
\newcommand{\rag}{\rangle}

\newcommand{\la}{\lambda}

\newcommand{\D}{\displaystyle}

\newcommand{\met}{\qquad \mbox{et} \qquad}

\newcommand{\saut}{\vspace{2ex}}
\newcommand{\vol}{\mathrm{Vol}}

\title{ An elementary solution to the Busemann-Petty problem}
\author{F. Barthe, M. Fradelizi and B. Maurey}
\date{}
\maketitle

\begin{abstract} A unified analytic solution to the Busemann-Petty problem
was recently found by Gardner, Koldobsky and Schlumprecht.
 We give an elementary proof of their
formulas for the inverse Radon transform of the radial function $\rho_K$ of 
an origin-symmetric star body $K$.
\end{abstract} 

Let $K$ and $L$ be two symmetric convex bodies 
in $\R^n$ such that for every hyperplane $H$ through the origin
\[ \vol_{n-1}(K \cap H) \le  \vol_{n-1}(L \cap H) \, ;\]
does it follow that
$  \vol_{n}(K) \le  \vol_{n}(L)\ ?$
The answer to this question of Busemann and Petty \cite{BP} is negative
for $n\ge 5$ (Gardner \cite{G1}, Papadimitrakis \cite{P}) and positive
for smaller dimensions (Gardner \cite{G2} for $n=3$, Zhang \cite{Z2}
for $n=4$). A unified solution was recently 
provided by Gardner, Koldobsky and Schlumprecht in \cite{GKS}, using
Fourier transform. We give an elementary proof of their
formulas for the inverse Radon transform of the radial function $\rho_K$.
\medskip

On $\R^n$, we denote the scalar product by 
$\lag \cdot,\cdot \rag$ and the Euclidean norm by $|\cdot|$\ .
We write $B^n$ for the unit ball and $S^{n-1}$ for the unit
sphere, and $v_{n}$, $s_{n-1}$ denote their respective volumes.
If $K\subset \R^n$ is a star body, its radial function 
$\rho_K$ is defined for every $x\in S^{n-1}$ by
\[ \rho_K(x)=\sup \left\{ \la >0 \, ;\; \la x \in K \right\}.\]
The connection between the Busemann-Petty problem and the spherical Radon
transform $R$ is due to Lutvak \cite{L}. Recall that $R$ acts on the 
space of continuous functions on  $S^{n-1}$ by setting
\[ Rf(\xi)=\int_{S^{n-1}\cap \xi^\bot} f(u) \, d\sigma_{n-2}(u) \]
for every $\xi \in S^{n-1}$; here $\sigma_{n-2}$ is the Haar measure 
of total mass $s_{n-2}$ on principal $n-2$ spheres. It follows from 
Lutvak \cite{L}, Zhang \cite{Z1}, that the Busemann-Petty problem has
a positive answer in $\R^n$ if and only if every  symmetric convex 
body $K$ in $\R^n$, with positive curvature and $C^\infty$ radial function, 
is such that $R^{-1}\rho_K$ is a non-negative function. In  
\cite{GKS}, the authors express $R^{-1}\rho_K$ in terms of 
\[ A_\xi(t)=\vol_{n-1}(K\cap(t\xi + \xi^\bot)), \quad \xi \in
S^{n-1}\]
as follows:
\medskip

{\bf \noindent Theorem } {\em Let $n\ge 3$.
Let $K\subset \R^n$ be an origin-symmetric star 
body, with $C^\infty$ radial function $\rho_K$.

If $n$ is even, then
\[ (-1)^{\frac{n-2}2}2^n \pi^{n-2} \rho_K 
        =R\left(\xi \mapsto A_\xi^{(n-2)} (0)\right).\]

If $n$ is odd, then
\[ \frac{(-1)^\frac{n-1}{2} (2\pi )^{n-1}}{(n-2)!}\rho_K
        = R\left(\xi \mapsto \int\limits_0^\infty t^{-n+1}\Bigl(
        A_\xi(t)-\sum_{k=0}^{\frac{n-3}2}  A^{(2k)}_\xi(0)\frac{t^{2k}}{(2k)!}
        \Bigr) dt \right).\]
}

{\bf \noindent Remark.} Let us recall why this solves the case $n = 4$
of the Busemann-Petty problem (\cite{Z2}, \cite{GKS}). If $n = 4$, then
$R^{-1}\rho_K(\xi) = - A''_\xi(0)/16\pi^2$. If $K$ is convex and
symmetric, the latter is non-negative (by Brunn-Minkowski,
the largest hyperplane section orthogonal
to $\xi$ is indeed the one through the origin).

\medskip
{\bf \noindent Proof.} We first compute the Radon transform of
$\xi \rightarrow A_\xi(t)$, for any given $t\ge 0$. 
Let $e\in S^{n-1}$ and set 
$f(t) := R(\xi\mapsto A_\xi(t))(e)$.
We identify $e^\bot$ and $\R^{n-1}$, and for $y\in \R^{n-1}$, 
we set $\phi(y) = \vol_1(K\cap (y + \R e))$. Then 
\beq
f(t)    & = & \int_{S_{n-1}\cap e^\bot}
        \int_{x \in \R^n,\,\lag x,\xi\rag=t} \mathbf{1}_K
        (x) \, d^{n-1}(x) \, d\sigma_{n-2}(\xi) \\
        & = & \int_{S_{n-1} \cap e^\bot}\int_{y\in e^\bot,\, \lag
        y,\xi\rag=t} \phi(y)\, d^{n-2}(y) \, d\sigma_{n-2}(\xi).
\eeq
Considered as a function of $g$, the quantity
\[ \int_{S_{n-1} \cap e^\bot}\int_{y\in e^\bot,\, \lag
        y,\xi\rag=t} g(y)\, d^{n-2}(y) \, d\sigma_{n-2}(\xi)\]
(where $g$ is defined on $e^\perp \simeq \R^{n-1}$)
is linear, continuous and rotation invariant. Hence 
there exists a measure $\mu_t$ on $\R^+$ such that for all 
$g$ the previous expression is equal to 
\[  \int_{\R^+} \left(\int_{S^{n-2}} g(ru)\,
d\sigma_{n-2}(u) \right) \, d\mu_t (r).\]
Applying the definition of $\mu_t$ with the function $g
=\mathbf{1}_{rB^{n-1}}$ yields 
\beq
s_{n-2}\ \mu_t([0,r])&=&\int_{S^{n-2}}\int_{\lag
        y,\xi\rag=t} \mathbf{1}_{rB^{n-1}}(y)\, d^{n-2}(y) \,
        d\sigma_{n-2}(\xi) \\
        &=& s_{n-2}\ v_{n-2}\  \mathbf{1}_{\{t\le r\}}
        (r^2-t^2)^{\frac{n-2}2}.
\eeq
Consequently, $\D d\mu_t(r)=s_{n-3}\
        r(r^2-t^2)^{\frac{n-4}2}\mathbf{1}_{\{t\le r\}}
        \, dr.$
Thus we have proved that
\[ f(t)= s_{n-3}\int_t^\infty r(r^2-t^2)^{\frac{n-4}2}\Phi(r)\,
dr,\]
where $\Phi$ is defined on $\R$ by
\[ \Phi(x) = \int_{S^{n-2}} \phi(xu)\, d\sigma_{n-2}(u).\]
Notice that $\Phi$ is even, compactly supported and $C^\infty$ in some 
neighborhood of the origin.  Our aim now is to relate $f(t)$ and
$\Phi(0) = 2\rho_K(e)s_{n-2}$. The case $n = 4$
is very simple: $f(t) = 2\pi\int_t^\infty r\Phi(r)dr$, 
hence $f''(0) = -2\pi \Phi(0) = -16\pi^2\rho_K(e)$.
By exchanging the order of the Radon transform
and the derivative, we conclude that $\rho_K$ is the Radon transform of  
$  \xi \mapsto -A''_\xi(0)/16\pi^2$.

\medskip
{\bf  If $n$ is even}:
\[
\frac{f(t)}{s_{n-3}}
        = \int_0^\infty r(r^2-t^2)^{\frac{n-4}2}\Phi(r)\,dr
         -t^{n-2}\int_0^1 u(u^2-1)^{\frac{n-4}2}\Phi(tu)\,du.
\]
The first term is a polynomial in $t$, of degree $n-4$ and 
$\Phi$ is $C^\infty$ in some neighborhood of $0$, thus
\[  f^{(n-2)}(0)= -s_{n-3}(n-2)!\int_0^1 u(u^2-1)^{\frac{n-4}2}\Phi(0)\,du
= (-1)^\frac{n-2}2 2^n \pi^{n-2} \rho_K(e).
\]
We conclude by exchanging the order of the Radon transform
and the derivative.

\medskip
{\bf If $n$ is odd}: the basic principle is still very simple, but
the technical details are slightly unpleasant. We shall begin by writing 
the proof as if $\Phi$ were $C^\infty$ on $\R\,$; but this 
is not true, because there are points of $e^\perp$ where 
our initial function $\phi$ is not differentiable, for example
the points of the boundary of the projection of $K$ on $e^\perp$; 
we shall indicate afterwards the standard approximation argument that 
fixes this difficulty. Integrating by parts, we get
\[F(t):=-\frac{n-2}{s_{n-3}}f(t)= \int_t^\infty(r^2-t^2)^\frac{n-2}2
        \Phi'(r)\, dr.\]
For $k\ge0$, let 
$a_k=(-1)^k {\frac{n-2}2 \choose k}
=\frac{(-1)^k}{k!}\prod_{j=0}^{k-1}(\frac{n-2}{2}-j)$. Notice that
$\sum |a_k|<\infty$. Let
\[ P(t)=\sum_{k=0}^{\frac{n-3}2} a_k t^{2k} \int_0^\infty r^{n-2-2k} 
        \Phi'(r) \, dr.\]
Then the quantity $\D \frac{F(t)-P(t)}{t^{n-1}}$ is equal to 
\[\int\limits_t^\infty \Bigl(  \sum_{k=\frac{n-1}2}^\infty a_k 
        \bigl(t^{-1}r\bigr)^{n-2-2k}\Bigr) \Phi'(r) \, \frac{dr}t
 -\int\limits_0^t \Bigl(  \sum_{k=0}^{\frac{n-3}2} a_k 
        \bigl(t^{-1}r\bigr)^{n-2-2k}\Bigr) \Phi'(r) \,
\frac{dr}t\]
\[=\int\limits_1^\infty \Bigl(  \sum_{k=\frac{n-1}2}^\infty a_k 
        u^{n-2-2k}\Bigr) \Phi'(tu) \, du 
         - \int\limits_0^1 \Bigl(  \sum_{k=0}^{\frac{n-3}2} a_k 
        u^{n-2-2k}\Bigr) \Phi'(tu) \, du.
\]
 By Fubini's theorem and since $\int_0^\infty \Phi'(tu)\,
dt=-\Phi(0)/u\,$, we get
\[\int_0^\infty  \frac{F(t)-P(t)}{t^{n-1}}\, dt =\Phi(0) \Bigl(
\sum_{k=0}^\infty  \frac{a_k}{n-2-2k} \Bigr)= c_n \rho_K(e),\]
which is finite. Thus, $P$ is the Taylor polynomial of $F$
of order $n-3$ at zero, and the above integral represents the action of 
the distribution $t_+^{-n+1}$ on $F$. We obtain therefore
\[ \lag t_+^{-n+1}, \, R(\xi \rightarrow A_\xi(t))(e) \rag =
- c_n \, {s_{n-3}\over n-2} \rho_K(e).
\]
A soft manner to compute $c_n$ is to replace $\Phi$ by 
$G(x) = \mathrm{e}^{-x^2}$ in the previous computation. 
Once again, we end the proof by exchanging the order 
in which the Radon transform and the distribution $t_+^{-n+1}$ 
act (we shall give some explanation about this at the end). 
\smallskip
 
 We now explain how to deal with the fact that $\Phi$ is
not $C^\infty$ everywhere. To every continuous and even function 
$\Phi_1$ on $\R$, which is $C^\infty$ in a neighborhood 
of $0$ and supported on a fixed interval 
$[-R, R]$ containing the support of $\Phi$, we associate the 
even function $F_1$ on $\R$ defined for $t \ge 0$ by
\[F_1(t) := -(n - 2) \int_t^\infty r(r^2 - t^2)^\frac{n-4}2
        \Phi_1(r) \, dr.\]
Let $Q(u)$ be the Taylor polynomial of degree $n - 3$ for
$(1 - u^2)^{(n-4)/2}$ at the origin, and let
$P_1(t) := -(n - 2) \int_0^\infty r^{n-3} Q(t/r) \Phi_1(r)\, dr$
(of course, $F_1 = F$ and $P_1 = P$ when $\Phi_1 = \Phi$). One can 
get easily the following estimates (where $C(n, R)$ or $C(a,n, R)$ 
denote constants depending only upon $n, R$ or $a, n, R$): 
\smallskip

 -- first, 
$\|F_1\|_\infty \le R^{n-2} \, \|\Phi_1\|_\infty$; 
\smallskip

 -- for every $t$, we have
$|P_1(t)| \le C(n,R) \, (1 + |t|^{n-3}) \|\Phi_1\|_\infty$; 
\smallskip

 -- finally, when $\Phi_1$ vanishes on some 
neighborhood $(-a, a)$ of $0$, one can see that
$|F_1(t) - P_1(t)| \le C(a, n, R) \, t^{n-1} \|\Phi_1\|_\infty$ 
for $0 \le t \le 1$. 
\smallskip

 These three estimates imply that the integral
$\int_0^\infty t^{-n+1} (F_1(t) - P_1(t)) \, dt$ converges 
to $\int_0^\infty t^{-n+1} (F(t) - P(t)) \, dt$ when we 
let $\Phi_1$, equal to $\Phi$ on a fixed interval $[-a, a]$ and 
supported on $[-R, R]$, tend uniformly to $\Phi$.
\medskip

 Let us turn finally to the interchange of the actions of the
Radon transform and the distribution $t_+^{-n+1}$ on the function 
$(\xi, t) \rightarrow A_\xi(t)$. It follows from our hypothesis that 
this function is $C^\infty$ on $S^{n-1} \times (-a, a)$ for some $a > 0$. 
Let us assume $n = 5$ for example. Since $K$ is symmetric, we may write
\[
A_\xi(t) = f_0(\xi) + t^2 f_2(\xi) + t^4 g(\xi, t)
\]
where $f_0$, $f_2$ and $g$ are continuous and bounded on $S^{n-1}$ and 
$S^{n-1} \times \R$ respectively. Since $A_\xi$ 
vanishes for $|t| > R$, we have 
$g(\xi, t) = - t^{-4} f_0(\xi) - t^{-2} f_2(\xi)$ for $t > R$, and 
\[
\lag A_\xi, t_+^{-4} \rag = 
\int_0^R g(\xi, t) \, dt - {R^{-3}\over 3} f_0(\xi) - R^{-1} f_2(\xi),
\]
which shows that the interversion with the integral over
$\xi\in S^{n-1}$ causes no trouble.
\qed

\bigskip

\bigskip




\hfill
\parbox{10cm}{\begin{center}
Equipe d'Analyse et Math\'ematiques Appliqu\'ees \\
Universit\'e de Marne la Vall\'ee \\
Boulevard Descartes, Cit\'e Descartes, Champs sur Marne\\
77454 Marne la Vall\'ee Cedex 2, FRANCE\\

\vspace{5mm}
\parbox{4.5cm}{
barthe@math.univ-mlv.fr\\
fradeliz@math.univ-mlv.fr\\
maurey@math.univ-mlv.fr}
\end{center}
}
\end{document}